\magnification=1200

\catcode`\ =\active \def { }

\hsize=11.25cm    
\vsize=18cm       
\parindent=12pt   \parskip=5pt     

\hoffset=.5cm   
\voffset=.8cm   

\pretolerance=500 \tolerance=1000  \brokenpenalty=5000

\catcode`\@=11

\font\eightrm=cmr8         \font\eighti=cmmi8
\font\eightsy=cmsy8        \font\eightbf=cmbx8
\font\eighttt=cmtt8        \font\eightit=cmti8
\font\eightsl=cmsl8        \font\sixrm=cmr6
\font\sixi=cmmi6           \font\sixsy=cmsy6
\font\sixbf=cmbx6

\font\tengoth=eufm10 
\font\eightgoth=eufm8  
\font\sevengoth=eufm7      
\font\sixgoth=eufm6        \font\fivegoth=eufm5

\skewchar\eighti='177 \skewchar\sixi='177
\skewchar\eightsy='60 \skewchar\sixsy='60

\newfam\gothfam           \newfam\bboardfam

\def\tenpoint{
  \textfont0=\tenrm \scriptfont0=\sevenrm \scriptscriptfont0=\fiverm
  \def\rm{\fam\z@\tenrm}
  \textfont1=\teni  \scriptfont1=\seveni  \scriptscriptfont1=\fivei
  \def\oldstyle{\fam\@ne\teni}\let\old=\oldstyle
  \textfont2=\tensy \scriptfont2=\sevensy \scriptscriptfont2=\fivesy
  \textfont\gothfam=\tengoth \scriptfont\gothfam=\sevengoth
  \scriptscriptfont\gothfam=\fivegoth
  \def\goth{\fam\gothfam\tengoth}
  
  \textfont\itfam=\tenit
  \def\it{\fam\itfam\tenit}
  \textfont\slfam=\tensl
  \def\sl{\fam\slfam\tensl}
  \textfont\bffam=\tenbf \scriptfont\bffam=\sevenbf
  \scriptscriptfont\bffam=\fivebf
  \def\bf{\fam\bffam\tenbf}
  \textfont\ttfam=\tentt
  \def\tt{\fam\ttfam\tentt}
  \abovedisplayskip=12pt plus 3pt minus 9pt
  \belowdisplayskip=\abovedisplayskip
  \abovedisplayshortskip=0pt plus 3pt
  \belowdisplayshortskip=4pt plus 3pt 
  \smallskipamount=3pt plus 1pt minus 1pt
  \medskipamount=6pt plus 2pt minus 2pt
  \bigskipamount=12pt plus 4pt minus 4pt
  \normalbaselineskip=12pt
  \setbox\strutbox=\hbox{\vrule height8.5pt depth3.5pt width0pt}
  \let\bigf@nt=\tenrm       \let\smallf@nt=\sevenrm
  \normalbaselines\rm}

\def\eightpoint{
  \textfont0=\eightrm \scriptfont0=\sixrm \scriptscriptfont0=\fiverm
  \def\rm{\fam\z@\eightrm}
  \textfont1=\eighti  \scriptfont1=\sixi  \scriptscriptfont1=\fivei
  \def\oldstyle{\fam\@ne\eighti}\let\old=\oldstyle
  \textfont2=\eightsy \scriptfont2=\sixsy \scriptscriptfont2=\fivesy
  \textfont\gothfam=\eightgoth \scriptfont\gothfam=\sixgoth
  \scriptscriptfont\gothfam=\fivegoth
  \def\goth{\fam\gothfam\eightgoth}
  
  \textfont\itfam=\eightit
  \def\it{\fam\itfam\eightit}
  \textfont\slfam=\eightsl
  \def\sl{\fam\slfam\eightsl}
  \textfont\bffam=\eightbf \scriptfont\bffam=\sixbf
  \scriptscriptfont\bffam=\fivebf
  \def\bf{\fam\bffam\eightbf}
  \textfont\ttfam=\eighttt
  \def\tt{\fam\ttfam\eighttt}
  \abovedisplayskip=9pt plus 3pt minus 9pt
  \belowdisplayskip=\abovedisplayskip
  \abovedisplayshortskip=0pt plus 3pt
  \belowdisplayshortskip=3pt plus 3pt 
  \smallskipamount=2pt plus 1pt minus 1pt
  \medskipamount=4pt plus 2pt minus 1pt
  \bigskipamount=9pt plus 3pt minus 3pt
  \normalbaselineskip=9pt
  \setbox\strutbox=\hbox{\vrule height7pt depth2pt width0pt}
  \let\bigf@nt=\eightrm     \let\smallf@nt=\sixrm
  \normalbaselines\rm}

\tenpoint

\def\pc#1{\bigf@nt#1\smallf@nt}         \def\pd#1 {{\pc#1} }

\catcode`\;=\active
\def;{\relax\ifhmode\ifdim\lastskip>\z@\unskip\fi
\kern\fontdimen2  -1.2 \fontdimen3 \string;}

\catcode`\:=\active
\def:{\relax\ifhmode\ifdim\lastskip>\z@\unskip\fi\penalty\@M\ \fi\string:}

\catcode`\!=\active
\def!{\relax\ifhmode\ifdim\lastskip>\z@
\unskip\fi\kern\fontdimen2  -1.1 \fontdimen3 \string!}

\catcode`\?=\active
\def?{\relax\ifhmode\ifdim\lastskip>\z@
\unskip\fi\kern\fontdimen2  -1.1 \fontdimen3 \string?}

\catcode`\«=\active 
\def«{\raise.4ex\hbox{%
 $\scriptscriptstyle\langle\!\langle$}}

\catcode`\»=\active 
\def»{\raise.4ex\hbox{%
 $\scriptscriptstyle\rangle\!\rangle$}}

\frenchspacing

\def\raggedbottom{\topskip 10pt plus 36pt\r@ggedbottomtrue}

\def\pointir{\unskip . --- \ignorespaces}

\def\Medbreak{\vskip-\lastskip\medbreak}

\long\def\th#1 #2\enonce#3\endth{
   \Medbreak\noindent
   {\pc#1} {#2\unskip}\pointir{\it #3}\smallskip}

\def\decale#1{\smallbreak\hskip 28pt\llap{#1}\kern 5pt}
\def\decaledecale#1{\smallbreak\hskip 34pt\llap{#1}\kern 5pt}
\def\puce{\smallbreak\hskip 6pt{$\scriptstyle\bullet$}\kern 5pt}

\def\eqalign#1{\null\,\vcenter{\openup\jot\m@th\ialign{
\strut\hfil$\displaystyle{##}$&$\displaystyle{{}##}$\hfil
&&\quad\strut\hfil$\displaystyle{##}$&$\displaystyle{{}##}$\hfil
\crcr#1\crcr}}\,}

\catcode`\@=12

\showboxbreadth=-1  \showboxdepth=-1

\newcount\numerodesection \numerodesection=1
\def\section#1{\bigbreak
 {\bf\number\numerodesection.\ \ #1}\nobreak\medskip
 \advance\numerodesection by1}

\mathcode`A="7041 \mathcode`B="7042 \mathcode`C="7043 \mathcode`D="7044
\mathcode`E="7045 \mathcode`F="7046 \mathcode`G="7047 \mathcode`H="7048
\mathcode`I="7049 \mathcode`J="704A \mathcode`K="704B \mathcode`L="704C
\mathcode`M="704D \mathcode`N="704E \mathcode`O="704F \mathcode`P="7050
\mathcode`Q="7051 \mathcode`R="7052 \mathcode`S="7053 \mathcode`T="7054
\mathcode`U="7055 \mathcode`V="7056 \mathcode`W="7057 \mathcode`X="7058
\mathcode`Y="7059 \mathcode`Z="705A


\catcode`\À=\active \defÀ{\`A}    \catcode`\à=\active \defà{\`a} 
\catcode`\Â=\active \defÂ{\^A}    \catcode`\â=\active \defâ{\^a} 
\catcode`\Æ=\active \defÆ{\AE}    \catcode`\æ=\active \defæ{\ae}
\catcode`\Ç=\active \defÇ{\c C}   \catcode`\ç=\active \defç{\c c}
\catcode`\È=\active \defÈ{\`E}    \catcode`\è=\active \defè{\`e} 
\catcode`\É=\active \defÉ{\'E}    \catcode`\é=\active \defé{\'e} 
\catcode`\Ê=\active \defÊ{\^E}    \catcode`\ê=\active \defê{\^e} 
\catcode`\Ë=\active \defË{\"E}    \catcode`\ë=\active \defë{\"e} 
\catcode`\Î=\active \defÎ{\^I}    \catcode`\î=\active \defî{\^\i}
\catcode`\Ï=\active \defÏ{\"I}    \catcode`\ï=\active \defï{\"\i}
\catcode`\Ô=\active \defÔ{\^O}    \catcode`\ô=\active \defô{\^o} 
\catcode`\Ù=\active \defÙ{\`U}    \catcode`\ù=\active \defù{\`u} 
\catcode`\Û=\active \defÛ{\^U}    \catcode`\û=\active \defû{\^u} 
\catcode`\Ü=\active \defÜ{\"U}    \catcode`\ü=\active \defü{\"u}

\def\hfl#1#2#3{\smash{\mathop{\hbox to#3{\rightarrowfill}}\limits
^{\textstyle#1}_{\textstyle#2}}}

\def\Q{{\bf Q}}

\def\R{{\bf R}}
\def\C{{\bf C}}
\def\N{{\bf N}}
\def\M{{\bf M}}
\def\H{{\bf H}}
\def\S{{\bf S}}

\def\Z{{\bf Z}}

\def\F{{\bf F}}

\def\Out{\mathop{\rm Out}\nolimits}
\def\Id{\mathop{\rm Id}\nolimits}

\def\End{\mathop{\rm End}\nolimits}

\def\Gal{\mathop{\rm Gal}\nolimits}

\def\zero{\{0\}}
\def\one{\{1\}}

\def\to{\rightarrow}

\centerline{\bf Numbers and periods}
\medskip
\centerline{\it Notes for a talk at the H-C.R.I., 8 February 2006}
\medskip
\centerline{Chandan Singh Dalawat}
\bigskip
We review the construction, starting with the monoid $\N$ of natural
numbers, of the ring $\Z$ of integers, of the field $\Q$ of rational
numbers, of the locally compact fields $\R$ and $\Q_p$ of real and
$p$-adic numbers ($p$ being a prime number), and of the complete
closed fields $\C$ and $\C_p$.  We introduce Fontaine's rings of
periods and, finally, indicate the need for going beyond them for the
purposes of outer galoisian actions.
\bigskip

{\bf 1}. ({\it Numbers as sums}).  Mathematics starts with 

\th DEFINITION 1
\enonce
The monoid\/ $\N$ is the free monoid on one generator, denoted\/ $1$.
\endth

We grant that such a monoid exists.  The elements of $\N$ are $0$,
$1$, $1+1$, $1+1+1$, $\ldots$; they are to be added for example as
$$
(1+1)+(1+1+1)=1+1+1+1+1
$$
with the requirement that $n+0=n$ and $0+n=n$ for every $n\in\N$.
Given any monoid $M$ and an element $m\in M$, there is a unique
morphism of monoids $f:\N\to M$ such that $f(1)=m$.

\indent\indent{\eightpoint Summary : $(\N,1)$ is the initial object in
the category of pointed monoids.}

\medbreak
{\bf 2}.  ({\it Numbers as products}).  Let $\End(\N)$ be the additive
monoid of endomorphisms of $\N$ :
$$
(f+g)(n)=f(n)+g(n).
$$
The map $e:\End(\N)\to\N$, $f\mapsto f(1)$ is an isomorphism of
monoids~; the reciprocal is the morphism $1\mapsto(1\mapsto 1)$,
i.e.~the unique morphism which sends $1$ to $\Id_\N$.

An endomorphism $f$ is injective if and only if $f(1)\neq0$, bijective
if and only if $f(1)=1$.  As the composite of two injective maps is
injective, $e$ can be used to make the set $\N^*$ (the complement of
$0$ in $\N$) into a monoid, by transport of structure~; it is {\it
commutative}.

For $m,n\in\N^*$, we say that $m$ divides $n$, and write $m|n$, if the
endomorphism $1\mapsto n$ factors via the endomorphism $1\mapsto m$.


\th DEFINITION 2
\enonce
A number\/ $p\in\N^*$ is called prime if\/ $p\neq1$ and if the only
divisors of\/ $p$ are\/ $1$ and\/ $p$.
\endth
The set of prime numbers is denoted by $P$.  The first few primes are
$2$, $3$, $5$, $7$, $11$, $\ldots$  Multiples of $2$ are called even
numbers, the others are called odd.  In particular, every prime $\neq2$
is odd.

\th THEOREM 1
\enonce
The map\/ $\N^{(P)}\to\N^*$ which send $f$ to $\prod_{p\in P}p^{f(p)}$
is an isomoprhism, i.e.~$\N^*$ is the free commutative monoid on the
set\/ $P$.
\endth

In other words, given any commutative monoid $M$ and any map $f:P\to
M$, there is a unique morphism $\N^*\to M$ which extends $f$.
Moreover, $P$ is the only subset of $\N^*$ with this property.

\th THEOREM 2
\enonce
The monoid\/ $\N^*$ is not finitely generated~: $P$ is infinite.
\endth

{\it Proof\/(\/{\rm Euclid}) :\/} Suppose that $P$ were finite, and
consider the number $n=1+\prod_{p\in P}p$, which is $>1$.  If $n$ is
prime, we get a contradiction because $P$ does not contain $n$.  If
$n$ is not a prime, we also get a contradiction, because it does not
have any prime divisor : no $p\in P$ is a divisor of $n$.

\indent\indent{\eightpoint Summary : $\N^*$ is the monoid of injective
endomorphisms of the monoid $\N$~; it is free commutative~; the
elements of its basis are called prime numbers.}

\medbreak
{\bf 3}. ({\it Numbers as differences}).  If $f:\N\to M$ is a morphism
of monoids such that $f(1)$ has an (additive) inverse, then $f(n)$ has
an inverse for every $n\in\N$.  Conversely, if $f(n)$ is invertible
for some $n\neq0$, then $f(1)$ is invertible.

The trouble is that there is no $n\in\N$ such that $n+1=0$, i.e.~the
difference $0-1$ does not exist.  Nor does the difference $8-9$, which
is essentially the same.  The way to make them exist is to declare
these differences $m-n$ ($m,n\in\N$) themselves as being new numbers
and to identify any two differences which are essentially the same.

Formally, consider the equivalence relation $\sim$ on the monoid
$\N\times\N$ of all pairs, defined as
$$
(m,n)\sim(m',n')\quad\Leftrightarrow\quad m+n'=m'+n.
$$
Then the class of $(0,1)$ is the same as the class of $(8,9)$, and so
on~; it is denoted $-1$.  The pairs $(1,0)$, $(9,8)$, etc., have the
same class, which is identified with $1$.  Thus the equivalence
classes consist of
$$
\ldots,\; -3,\; -2,\; -1,\; 0,\; 1\;, 2,\; 3,\ldots
$$
and we have obtained the group $\Z$ of integers as the group of
differences of the monoid $\N$

\th THEOREM 3
\enonce
The group\/ $\Z$ is free on one generator~; $1$ and $-1$
are the only generators.  
\endth

\indent\indent{\eightpoint Summary : $\N\to\Z$ is the group of
differences of $N$, i.e.~the initial object in the category of maps of
the monoid $\N$ into groups.}

\medbreak
{\bf 4}.  ({\it Numbers as products, bis}).  Consider the ring
$\End(\Z)$ of endomorphisms of the group $\Z$.  The map
$\End(\Z)\to\Z$, $f\mapsto f(1)$ is an isomorphism of additive groups,
allowing us to make the commutative group $\Z$ into a ring~; extending
the multiplication from $\N$ to $\Z$.  This ring is commutative and
integral, i.e.,~if $xy=0$ ($x,y\in\Z$), then $x=0$ or $y=0$.
Moreover, the sum and the product of two positive numbers is again
positive.

\th THEOREM 4
\enonce
The ring\/ $\Z$ is principal, i.e., every ideal can be generated by
one element.
\endth

Given any ring $A$, there is a unique homomorphism of rings $\Z\to A$,
a fact which is expressed by saying that the ring $\Z$ is the initial
object in the category of rings.

\th THEOREM 5
\enonce
For every prime~$p$, the ideal\/ $p\Z$ is prime, indeed maximal.
Conversely, every maximal ideal of\/ $\Z$ is generated by a prime
number.
\endth

Notice that $0$ is the only prime ideal of $\Z$ which is not maximal.
For every prime~$p$, the field $\Z/p\Z$ is denoted $\F_p$.

\indent\indent{\eightpoint Summary : The ring $\Z$ is the ring of
endomorphisms of the (commutative) group $\Z$.  It is principal and
the maximal ideal correspond to prime numbers.}

\medbreak
{\bf 5}.  ({\it Numbers as fractions}).  We notice that fractions do
not exist in $\Z$: there is no $x$ such that $2x=1$.  In order to
remedy this, we proceed as with difference : on the set $\Z\times\Z^*$
of pairs $(m,n)$ ($n\neq0$), we define an equivalence relation $\sim$
as
$$
(m,n)\sim(m',n')\quad\Leftrightarrow\quad mn'=m'n.
$$
The equivalence class of $(m,n)$ is written $m/n$~; they are to be
added and multiplied as
$$
{m\over n}+{m'\over n'}={mn'+m'n\over nn'},\qquad 
{m\over n}\cdot{m'\over n'}={mm'\over nn'}.
$$
These classes form a field $\Q$, called the field of fractions of the
(integral) ring $\Z$.

Consider the multiplicative group $\Q^\times$ of the field $\Q$, and
its subgroup $\Z^\times$, which happens to be the torsion subgroup.
One identifies $P$ with its image in $\Q^\times\!/\Z^\times$.

\th THEOREM 6
\enonce
The map\/ $\Z^{(P)}\to\Q^\times\!/\Z^\times$ which sends $f$ to
$\prod_{p\in P}p^{f(p)}$ is an isomoprhism, i.e.
$\Q^\times\!/\Z^\times$ is the free $\Z$-module on the set\/ $P$.
\endth
In other words, $\Q^\times\!/\Z^\times$ is the group of
``differences'' of the monoid $\N^*$.

\th COROLLARY
\enonce
The map\/ $\{1,-1\}\times\Z^{(P)}\to\Q^\times$ which sends\/
$(\varepsilon,f)$ to $\varepsilon\prod_{p\in P}p^{f(p)}$ is an
isomoprhism of\/ $\Z$-modules.
\endth

Sometimes it is more convenient to use the isomorphism which sends
$-1$ to $-1$ and $e_p$ to $p^*$ for every prime $p$, where
$(e_p)_{p\in P}$ is the canonical basis of $\Z^{(P)}$ and
$$
2^*=2,\qquad p^*=p^{p-1\over2}\quad(p\neq2).
$$

\indent\indent{\eightpoint Summary : $\Z\to\Q$ is the field of fractions
of $\Z$, i.e.~the initial object in the category of injective
homomorphisms of the ring $\Z$ into fields.
$\N^*\to\Q^\times\!/\Z^\times$ is the ``group of differences'' of
$\N^*$.}.

\medbreak
{\bf 6}.  ({\it Numbers as limits}).  We notice that there are
sequences $(x_n)_n$ of rational numbers which do not have a limit in
$\Q$, although $|x_m-x_n|$ can be made as small as we please by taking
$m,n$ large enough.  For example,
$$
x_n={1\over0!}+{1\over1!}+{1\over2!}+{1\over3!}+\cdots+{1\over n!}
$$
In order to make such limits exist, we consider the space of all
fundamental sequences, i.e.~sequences $(x_n)_n$ of rational numbers
such that $|x_n-x_m|$ can be made as small as we please by taking
$m,n$ large enough.  They can be added and multiplied, so they form a
ring.  Those fundamental sequences which tend to $0$, i.e.~such that
$|x_n|$ is as small as we please for $n$ large enough, form a maximal
ideal of this ring.  The quotient field $\R$ has a natural order,
hence a uniformity and a topology.
\th THEOREM 7
\enonce
The space\/ $\R$ is connected, locally connected, locally compact and
complete.
\endth

\medbreak
{\bf 7}.  ({\it Numbers as roots}).  The polynomial $X^2+1$ is
irreducible over $\R$ because $x^2+1>0$ for every $x\in\R$.  Adjoining
a root of this polynomial, we obtain the field $\C=\R[X]/(X^2+1)$ of
complex numbers.  Something miraculous happens~:
\th THEOREM 8
\enonce
The field\/ $\C$ of complex numbers is algebraically closed.
\endth
The proof rests on the fact that every real number $>0$ has a square
root, and that every real polynomial of odd degree has a root in $\R$.

The group of $\R$-automorphisms of $\C$ is generated by the involution
$i\mapsto-i$, where $i$ is the image of $X$ in $\C$.

The field $\C$ comes with an absolute value $|\phantom{x}|_\infty$
extending the one on $\R$, which can be used to give it a uniformity,
and hence a topology.  Another miracle happens~:

\th THEOREM 9
\enonce
The space\/ $\C$ is simply connected, locally compact and complete.
\endth

Let us mention the field $\H$ of quaternions, which is not
commutative.  It was discovered by Hamilton, who was trying to find an
even bigger number system than $\C$.  It can now be defined as the
unique ($4$-dimensional) $\R$-algebra $A$, other than $\M_2(\R)$, such
that $A\otimes_\R\C$ is $\C$-isomorphic to $\M_2(\C)$.  The absolute
value $|\phantom{x}|_\infty$ extends to $\H$, with respect to which it
is complete.
\th THEOREM 10
\enonce
The only fields complete with respect to an archimedean absolute value
are $\R$, $\C$ and $\H$.
\endth
In fact, these are the only locally compact connected fields.  Notice
that these three fields are different, as $\R$ is commutative but not
algebraically closed, $\C$ is commutative and algebraically closed,
whereas $\H$ is not commutative.

\indent\indent{\eightpoint Summary : $\Q\to\R$ is the completion of $\Q$
in the $\infty$-adic uniformity.  It is connected, locally connected,
locally compact and complete. $\R\to\C$ is a quadratic extension and
an algebraic closure, with a chosen $4^{\rm th}$ root of ~$1$.}

\medbreak
{\bf 8}.  ({\it Numbers as limits, bis\/}).  It was noticed at the end
of the XIX$^{\rm th}$ century that $|\phantom{x}|_\infty$ is not the
only absolute value on $\Q$.  For every prime $p$, there is an
absolute value $|\phantom{x}|_p$ for which the sequence
$$
x_n=1+p+p^2+p^3+\cdots+p^n
$$
is a fundamental sequence~!  It is the unique homomorphism
$\Q^\times\to\R^{\times\circ}$ (into the multiplicative group of real
numbers $>0$) such that $|p|_p=1/p$ and $|l|_p=1$ for every prime
$l\neq p$ (cf.~Corollay to Theorem~6).

\th THEOREM 11
\enonce
Leaving aside the trivial one, every absolute value on $\Q$ is
equivalent to $|\phantom{x}|_\infty$ or to $|\phantom{x}|_p$ for some
prime $p\in P$.
\endth

That none of these absolute values can be neglected is illustrated by
the following result (the ``product formula'') --- almost a tautology
--- in which $\bar P$ denoted the set $P$ together with an additional
element $\infty$~:

\th THEOREM 12
\enonce
For every $x\in\Q^\times$, one has\/ $|x|_v=1$ for almost all\/
$v\in\bar P$ and\/ $\prod_{v\in\bar P}|x|_v=1$.
\endth

Indeed, it is sufficient to verify this for $x=-1$ and $x=p$ ($p\in
P$), cf.~Corollary to Theorem~6.

There are fundamental sequences with respect to $|\phantom{x}|_p$
which do not converge in $\Q$, for example the sequence $x_n=2^{5^n}$
for $p=5$, whose limit would be a primitive $4^{\rm th}$ root of~$1$.

As in the construction of real numbers, we can consider the ring of
all fundamental sequences~; those which converge to $0$ form a maximal
ideal.  The quotient field $\Q_p$ comes with an absolute value
$|\phantom{x}|_p$ (extending the one on $\Q$), hence a metric.

\th THEOREM 13
\enonce
The space $\Q_p$ is locally compact, totally disconnected, and
complete.
\endth

The disc $|x|_p\le1$ is a compact subring $\Z_p$, with group of units
$\Z_p^\times$ the circle $|x|_p=1$.  Another simplification occurs in
the multiplicative group (cf.~Theorem~6).

\th THEOREM 14
\enonce
The valuation $v_p:\Q_p^\times\to\Z$ induces an isomorphism
$\Q_p^\times/\Z_p^\times\to\Z$. 
\endth

Choosing a section of $v_p$, i.e.~an element $e_p\in\Q^\times$ such
that $v_p(e_p)=1$ (for example $e_p=p$), we get an isomorphism
$\Z_p^\times\times\Z\to\Q_p^\times$.  For $p$ odd, the torsion
subgroup $W$ of $\Z_p^\times$ is cyclic of order~$p-1$~; the quotient
$\Z_p^\times/W$ is a commutative pro-$p$-group, i.e.~a $\Z_p$-module~;
it is free of rank~$1$~; the image of $1+e_p$ is a basis.  What
happens for $p=2$~?

What is the archimedean analogue of this theorem~?  The analogue of
the units would be the kernel of the absolute value, i.e.~the
$0$-sphere $\{1,-1\}$ in $\R$, the $1$-sphere $\S_1$ in $\C$ and the
$3$-sphere $\S_3$ in $\H$.

\th THEOREM 15
\enonce
The maps $(\varepsilon,x)\mapsto\varepsilon e^x$,
$\{1,-1\}\times\R\to\R^\times$, $\S_1\times\R\to\C^\times$ and\/
$\S_3\times\R\to\H^\times$ are isomorphisms of real-analytic groups.
\endth

In other words, the topological groups\/ $\R^\times\!/\{1,-1\}$,
$\C^\times\!/\S_1$ and\/ $\H^\times\!/\S_3$ are\/ $1$-dimensional
vector spaces over\/ $\R$, with canonical basis\/ $\bar e$~; they are
the same as the multiplicative group $\R^{\times\circ}$ of reals $>0$.

Further, there is a unique $\pi>0$ in $\R$ such that $2i\pi\Z$ is the
kernel of $e^{(\phantom{x})}:2i\pi\R\to\S_1$, i.e.~one has the exact
sequence
$$
\zero\to2i\pi\Z\to2i\pi\R\to\S_1\to\one.
$$
What is the analogue for $\S_3$~?  What is the liegebra ${\goth
s}_3$~?  What are the groups $\C^\times\!/\R^\times$,
$\H^\times\!/\R^\times$, $\H^\times\!/\{1,-1\}$, $\S_1/\{1,-1\}$,
$\S_3/\{1,-1\}$~?

In all fairness, we should give a similar description of the
multiplicative group of all fields complete with respect to the
$p$-adic absolute value, and not just $\Q_p$.  Let us leave that
aside.

\indent\indent{\eightpoint Summary : $\Q\to\Q_p$ is the completion of
$\Q$ in the $p$-adic uniformity~; it is locally compact and totally
disconnected.  The closed unit disc $\Z_p$ is a compact subring~; the
quotient $\Q_p^\times\!/\Z_p^\times$ is $\Z$.}

\medbreak
{\bf 9}. ({\it Numbers as roots and limits, bis\/}).  We saw that
there is a drastic reduction in the number of irreducible polynomials
when we completed $\Q$ with respect to $|\phantom{x}|_\infty$.  The
reduction is not so drastic when we complete with respect to
$|\phantom{x}|_p$~: there are irreducible polynomials of every degree,
but essentially only finitely many of them~:
\th THEOREM 16
\enonce 
The field $\Q_p$ has only finitely many extensions in each degree.
\endth
For example, the only quadratic extensions of $\Q_5$ are
$\Q_5(\sqrt{5})$, $\Q_5(\sqrt{2})$ and $\Q_5(\sqrt{10})$.

Central simple $\Q_p$-algebras are classified by the group $\Q/\Z$,
whereas central simple $\R$-algebras are classified by the group
${1\over2}\Z/\Z$.  There is a very precise local-to-global principle
in terms of these groups, of which a special case reads~:
\th THEOREM 17
\enonce
An\/ $n^2$-dimensional central simple\/ $\Q$-algebra\/ $A$ is
isomorphic to\/ $\M_2(\Q)$ if\/ $A\otimes_\Q\Q_p$ is isomorphic to\/
$\M_2(\Q_p)$ for every\/ $p\in P$~; $A\otimes_\Q\R$ is then
necessarily isomorphic to\/ $\M_2(\R)$.
\endth

Let $\overline\Q_p$ be an algebraic closure of $\Q_p$~; the group
$\Gal(\overline\Q_p|\Q_p)$, being pro-solvable, is simpler than
$\Gal(\overline\Q|\Q)$, but it is considerably more complicated than
the \hbox{order-$2$} group $\Gal(\C|\R)$.  The absolute value
$|\phantom{x}|_p$ extends uniquely to $\overline\Q_p$.  Here is a
minor embarassment~:
\th THEOREM 18
\enonce
The field\/ $\overline\Q_p$ is not complete (and not locally compact).
\endth

But we can complete $\overline\Q_p$ with respect to $|\phantom{x}|_p$
to obtain a field $\C_p$, to which $|\phantom{x}|_p$ extends uniquely,
with the same group of values.
\th THEOREM 19
\enonce
The field $\C_p$ is algebraically closed.
\endth
So we don't need to start all over again, taking an algebraic
closure, etc.

\indent\indent{\eightpoint Summary : $\overline\Q_p$ is an algebraic
closure of $\Q_p$.  $\C_p$ is the completion of $\overline\Q_p$ in the
$p$-adic uniformity; it is algebraically closed.}

\medbreak
{\bf 10}. ({\it Numbers as periods\/}).  Varieties have periods.  Let
us return briefly to the archimedean world of $\R$ and $\C$ to
illustrate this.

A considerable effort has been spent, following Grothendieck's
insights, in the search for a universal cohomology theory for
varieties over $\Q$.  It was noticed that the various known cohomology
theories share a number of properites, which led Grothendieck to
suspect that there must be a universal cohomology theory --- a theory
of motives ---, of which the known choosing theories are but various
different ``avatars''.  In other words, the motive of a variety should
carry every bit of cohomological information about the variety.

To be specific, let us take a smooth, projective, absolutely connected
curve $X$ of genus $g$ over $\R$.  Then $H^1(X(\C),\Z)$ is a free
$\Z$-module of rank $2g$, so all it knows about $X$ is its genus.
However, the $\C$-space $H^1(X(\C),\Z)\otimes\C$ has a natural
filtration which carries much more information about $X$, and indeed
enough information to recover the jacobian $J$ of $X$.  What I want to
emphasize is that you need to tensor with $\C$ to get this additional
structure.  

The miracle is that tensoring with $\C$ suffices for varieties over
$\R$.  Another way of saying this is that the field $\C$ contains the
periods of all smooth projective $\R$-varieties.  This is a miracle
because $\C$ was not designed to serve this purpose~; it was merely
designed to have a square root of $-1$.

\medbreak

{\bf 11}. ({\it Numbers as periods, bis\/}).  Varieties over $\Q_p$
have their periods too.  The major difference with the archimedean
world is that you have to go way beyond $\C_p$ to get all the periods.
Fontaine has pursued relentlessly this quest for the $p$-adic analogue
of $2i\pi$~; his solution is to construct several rings carrying many
structures which play the role of rings of periods for different kinds
of varieties~: those with smooth reduction, those with semistable
reduction, all (smooth proper) varieties.

The rings $B$ constructed by Fontaine are $\Q_p$-algebras which carry
--- among other structures --- a representation of
$G_p=\Gal(\overline\Q_p|\Q_p)$~; the invariants $B^{G_p}$ are a field.
A $p$-adic representation $E$ of $G_p$, i.e.~a finite-dimensions
vector $\Q_p$-space with a continuous liner action of $G_p$, is called
$B$-admissible if the dimension of $(B\otimes_{\Q_p}E)^{G_p}$ (as a
vector space over $B^{G_p}$) is the same as the dimension of $E$ (as a
vector space over $\Q_p$).  Notable examples of $p$-adic
representations are provided by the $p$-adic {\'e}tale cohomology of
varieties over $\Q_p$ --- the analogue of the singular cohomology of
varieties over $\R$ and $\C$.

The rings of Fontaine have proved their utility beyond any doubt.  A
theorem of Colmez gives an idea of their coherence as $p$ varies~: he
proves a product formula (cf.~Theorem~12) for the periods of abelian
varieties having complex multiplications, periods which may be
trascendental~!  

I have chosen to illustrate their utility by taking for $B$ the
crystalline ring of Fontaine~; $B$-admissible representations of $G_p$
are called crystalline.

\th THEOREM (Faltings)
\enonce
For a smooth projective variety\/ $X$ having smooth reduction over
$\Q_p$, the representation\/ $H^n(X_{\overline\Q_p},\Q_p)$ is
crystalline.
\endth

The converse need not hold.  Nevertheless, in the light of the
Serre-Tate criterion (see below), one has a right to
expect it to hold for abelian varieties.

\th THEOREM (Coleman-Iovita)
\enonce
Let $A$ be an abelian $\Q_p$-variety.  If the $p$-adic
representation\/ $H^1(A_{\overline\Q_p},\Q_p)$ is crystalline, then
$A$ has abelian reduction.
\endth

If we compare these two theorems with what goes on for $l$-adic
{\'e}tale cohomology, where $l$ is a prime $\neq p$, the similarity is
striking~:

\th THEOREM (Grothendieck)
\enonce
For a smooth projective variety\/ $X$ having smooth reduction over
$\Q_p$, the representation\/ $H^n(X_{\overline\Q_p},\Q_l)$ is
unramified.
\endth

\th THEOREM (Serre-Tate)
\enonce
Let $A$ be an abelian $\Q_p$-variety.  If the $l$-adic
representation\/ $H^1(A_{\overline\Q_p},\Q_l)$ is unramified, then $A$
has abelian reduction.
\endth

\indent\indent{\eightpoint Summary :  For abelian varities over a finite
extension of $\Q_p$, there is an $l$-adic and a $p$-adic criterion for
abelian reduction.}


\medbreak

{\bf 12}.  ({\it Anabelian periods\/}).  These criteria of Serre-Tate
and Coleman-Iovita fail for curves.  Indeed, it is not very difficult
to construct smooth proper absolutely connected curves $C$ over $\Q_p$
which do not have smooth reduction but whose jacobian has abelian
reduction.  For such a curve $C$, the $l$-adic {\'e}tale cohomology
$H^1(C_{\overline\Q_p},\Q_l)$ is unramified and the $p$-adic {\'e}tale
cohomology $H^1(C_{\overline\Q_p},\Q_p)$ is crystalline, because these
groups are the same for $C$ and its jacobian.

Is there a criterion for the smooth reduction of curves~?  A very
striking theorem was proved by Takayuki Oda, who gives an anabelian
$l$-adic criterion.  He replaces $H^1(C_{\overline\Q_p},\Q_l)$ by the
maximal pro-$l$-quotient $\pi_1^{(l)}(C_{\overline\Q_p})$ of the
fundamental group $\pi_1(C_{\overline\Q_p})$ of $C_{\overline\Q_p}$
(there is a choice of a base point involved~; it will turn out to be
immaterial).  The fundamental group carries a natural outer action of
the Galois group, i.e.~there is a natural homomorphism
$G_p\to\Out\pi_1^{(l)}(C_{\overline\Q_p})$.

\th THEOREM (Takayuki Oda)
\enonce
Let\/ $C$ be a smooth proper absolutely connected\/ $\Q_p$-curve of
genus\/ $>1$.  If the the outer action\/
$G_p\to\Out\pi_1^{(l)}(C_{\overline\Q_p})$ on the pro-$l$ $\pi_1$ is
unramified, then\/ $C$ has smooth reduction.
\endth
In fact, his theorem works over any finite extension of $\Q_p$.

\indent\indent{\eightpoint Summary : For curves, there is an anabelian $l$-adic
criterion for smooth reduction.}

\medbreak

{\bf 13}.  ({\it Anabelian periods, bis\/}).  Let me end with a
problem which I would very much like to solve.  Anyone who has read
this account cannot fail to ask himself~: Is there an anabelian
$p$-adic criterion for smooth reduction of curves over $\Q_p$~?  I
happened to mention this to Kazuya Kato, and he said that this problem
``goes straight to the heart''.  

\th PROBLEM
\enonce
Give a necessary and sufficient condition on the outer action
$G_p\to\Out\pi_1^{(p)}(C_{\overline\Q_p})$ for the (smooth,
projective, absolutely connected) $\Q_p$-curve $C$ (of genus $>1$) to
have smooth reduction.
\endth

In other words, one should complete the square~:
\smallskip

\setbox11=\vbox{\halign{\hfil#\hfil\cr
 Serre-Tate\cr
 \ \ (1968, ``unramified'')\ \ \cr}}

\setbox12=\vbox{\halign{\hfil#\hfil\cr
 Coleman-Iovita\cr
 \ \ (1999, ``crystalline'')\ \ \cr}}

\setbox21=\vbox{\halign{\hfil#\hfil\cr
 Takayuki Oda\cr
 \ \ (1995, ``unramified'')\ \ \cr}}

\setbox22=\vbox{\halign{\hfil#\hfil\cr
 ???\cr
 (???)\cr}}

\smallskip

\halign{
\hfil#\hfil\qquad& 
\hfil#\hfil& 
\hfil#\hfil\cr 
&$l\neq p$&$l=p$\cr
\noalign{\smallskip\hrule\smallskip}
abelian varieties&\box11&\box12\cr
\noalign{\smallskip\hrule\smallskip}
curves of genus $>1$&\box21&\box22\cr
\noalign{\smallskip\hrule}
}
\smallskip

The problem should be studied not just over $\Q_p$, but over every
finite extension thereof.

\indent\indent{\eightpoint Summary : Find an anabelian $p$-adic
criterion for smooth reduction of curves.}

{\bf 14} ({\it Semistable reduction}).  Notice that there is already a
$p$-adic criterion for semistable reduction of curves and abelian
varieties in terms of Fonaine's semistable ring.  In effect, Fontaine
has proved that if an abelian $K$-variety has semiabelian reduction,
$K$ being a finite extension of $\Q_p$, then the representation
$H^1(A_{\overline K},\Q_p)$ of $\Gal(\overline K|K)$ is semistable,
$\overline K$ being an algebraic closure of $K$.  Later, Breuil proved
the converse, at least for $p\neq2$.  Earlier, Deligne and Mumford had
proved that for a curve to have semistable reduction, it is necessary
and sufficient for its jacobian to have semiabelian reduction.

\bye